\documentclass[12pt]{article}
\usepackage{graphicx,amssymb,amsmath}
\usepackage{hyperref,url,breakurl}
\usepackage{verbatim}
\newcommand{\nocopyright}{
No Copyright\thanks{
The authors hereby waive all copyright
and related or neighboring rights to this work,
and dedicate it to the public domain.
This applies worldwide.
}}
\title{Frobenius's last proof}
\author{Peter G. Doyle}
\date{Version 1.0 dated 13 April 2019\\
\nocopyright
}
\newtheorem{theorem}{Theorem}
\newtheorem{prop}[theorem]{Proposition}

\newcommand{\proofstart}{{\noindent \bf Proof.\ }}
\newcommand{\qed}{\spadesuit}
\newcommand{\mathproofend}{\quad \qed}
\newcommand{\proofend}{$\quad \qed$}
\newcommand{\goesto}{\rightarrow}

\newcommand{\fig}[3]{
\begin{figure}
\includegraphics[width=370pt]{figures/#1.pdf}
\caption{#3}
\label{#2}
\end{figure}
}
\newcommand{\figsize}[4]{
\begin{figure}
\centerline{
\includegraphics[width=#1]{figures/#2.pdf}
}
\caption{#4}
\label{#3}
\end{figure}
}

\newcommand{\gversus}[2]{\left[#1\,|\,#2\right]}

\newcommand{\floor}[1]{\left \lfloor #1 \right \rfloor}

\newcommand{\gbinom}{\genfrac{[}{]}{0pt}{}}

\begin{document}

\maketitle

\begin{abstract}
Around about 1917,
Issai Schur
rediscovered the Rogers-Ramanujan identities,
and proved a system of polynomial identities that imply them.
Schur wrote that 
Georg Frobenius (his former advisor)
had shown him a simple, direct proof of these polynomial identities.
Schur did not see fit to reveal Frobenius's proof,
preferring his own rather complicated proof.
But it is easy enough to guess what this `simple, direct' proof must have been.
As Frobenius died in 1917, we may call this `Frobenius's last proof'.
\end{abstract}

\section{Introduction}

Around about 1917,
Issai Schur
\cite{schur:rr}
rediscovered the Rogers-Ramanujan identities,
and proved a system of polynomial identities that imply them.
Schur
\cite[p. 131]{schur:rr}
wrote that 
Georg Frobenius (his former advisor)
had shown him a simple, direct proof
(`einen enfachen direkten Beweis')
of these polynomial identities:
Schur did not see fit to reveal Frobenius's proof,
preferring his own rather complicated proof.
But it is easy enough to guess what this `simple, direct' proof must have been.
As Frobenius died in 1917, we may call this `Frobenius's last proof'.

\section{The alpine identities}

Told that Figure \ref{fig:gs1}
illustrates an identity involving binomial coefficients,
along with the hint that the specific colors red and green are significant,
many a sophomore will be able to guess the identity.

\fig{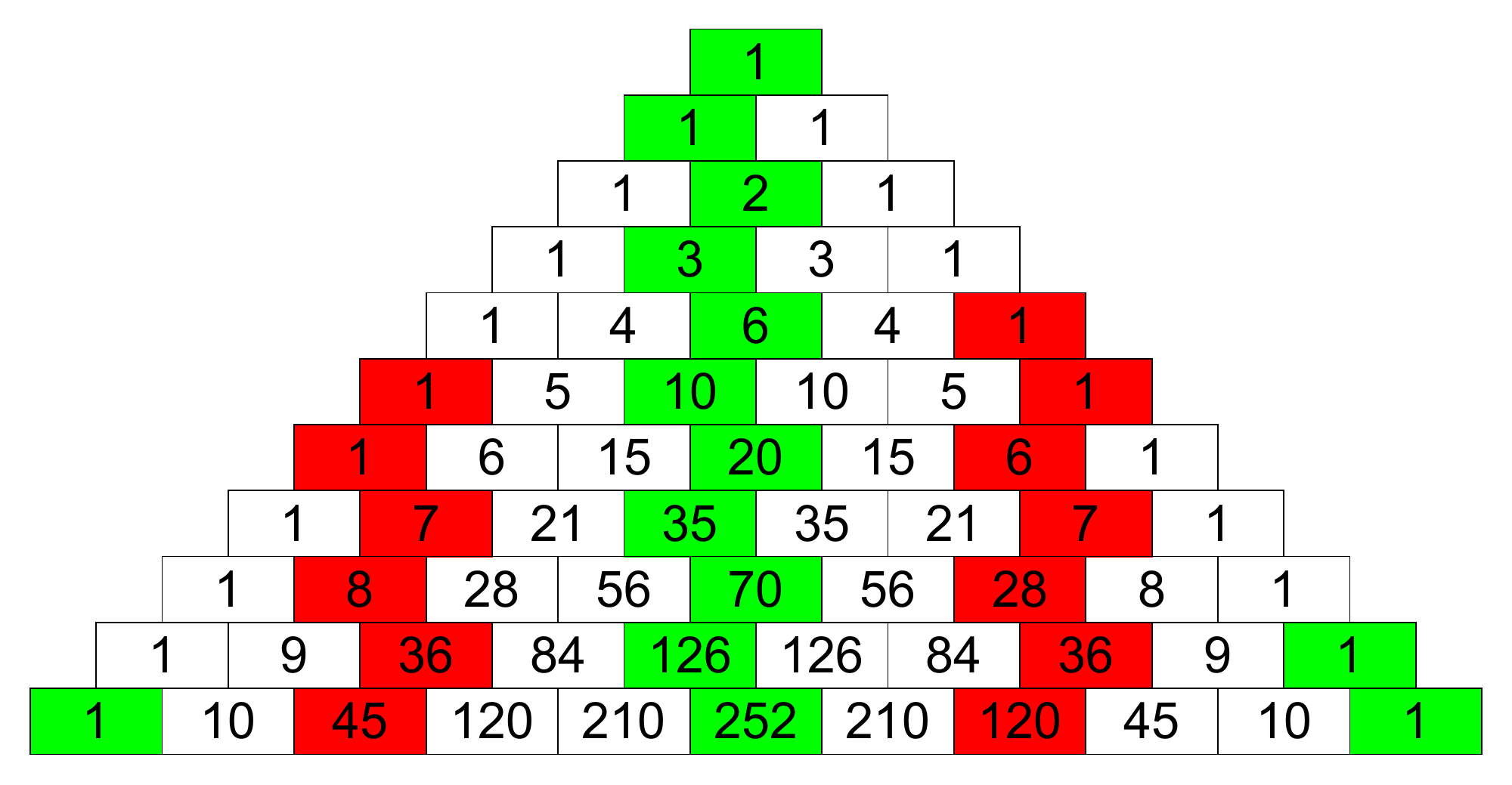}{fig:gs1}{First GS identity.}

What gives sophomores a leg up here
is that the Fibonacci sequence is apt to loom large in their imagination.
This steers them to the answer:
Totaling the numbers in each row,
taking the green numbers as positive and the red as negative,
yields the Fibonacci sequence.

Schur
\cite[p. 131]{schur:rr}
thanks Frobenius for showing him the following
way to formulate the identity.
(We're shifting the sequence one notch from what is in
Schur, to make $p(0)=p(1)=1$.)
\begin{prop}[First GS identity --- amateur version]
The sequence
\[
p(n)
=
\sum_{
-\floor{\frac{n}{5}}
\leq \lambda \leq
\floor{\frac{n+1}{5}}
}
(-1)^\lambda
\binom{n}{\floor{\frac{n+5\lambda}{2}}}
\]
satisfies the recurrence
\[
p(n) = p(n-1)+p(n-2)
\]
with
\[
p(-1)=0;\; p(0)=1
.
\]
\end{prop}
In the summation here, $\lambda$ may range over all integers;
the limits simply indicate the non-vanishing terms.

We're calling this the `GS identity',
short for `Giant Slalom identity',
because of the way the highlighted entries slalom their way
down Pascal's triangle.
This is the `amateur version' of the identity,
formulated for standard binomial coefficients
$\binom{n}{k}$.
Schur's original `pro version' involves
Gauss's polynomial analogs $\gbinom{n}{k}$ of the binomial coefficients.
We start with the amateur identity because
the pro identity is just a gussied-up version of the amateur identity.
The same goes for the proof.

And as for the proof,
many a sophomore
(though this time, not as many)
can find it.
It is illustrated in Figure \ref{fig:gs1proof}.
\fig{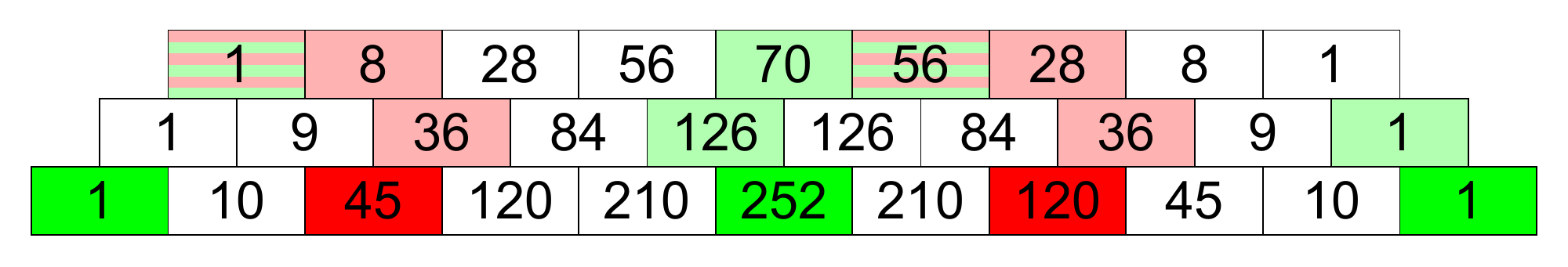}{fig:gs1proof}{Proof of the first GS identity.}
This proof is as simple and direct as one could wish.
We unhesitatingly identify this 
as `Frobenius's last proof'.

In symbols, the proof comes down to this:
\begin{eqnarray*}
\binom{n}{k-1}
&=&
\binom{n-1}{k-2} + \binom{n-1}{k-1}
\\&=&
\binom{n-1}{k-2}
+ \binom{n-2}{k-2} + \binom{n-2}{k-1}
,
\end{eqnarray*}
while
\begin{eqnarray*}
\binom{n}{k+1}
&=&
\binom{n-1}{k} + \binom{n-1}{k+1}
\\&=&
\binom{n-2}{k-1} + \binom{n-2}{k}
+ \binom{n-1}{k+1}
.
\end{eqnarray*}
Subtracting, we have 
\begin{eqnarray*}
\binom{n}{k-1} - \binom{n}{k+1}
&=&
\binom{n-1}{k-2}
- \binom{n-1}{k+1} +
\\&&
\binom{n-2}{k-2}
- \binom{n-2}{k}
,
\end{eqnarray*}
because the terms $\binom{n-2}{k-1}$ cancel.

As the name `First GS identity' indicates, there is a second GS identity,
illustrated in Figure \ref{fig:gs2}.
Here again the row totals give the Fibonacci sequence, this time
starting with $0,1$ instead of $1,1$.

\fig{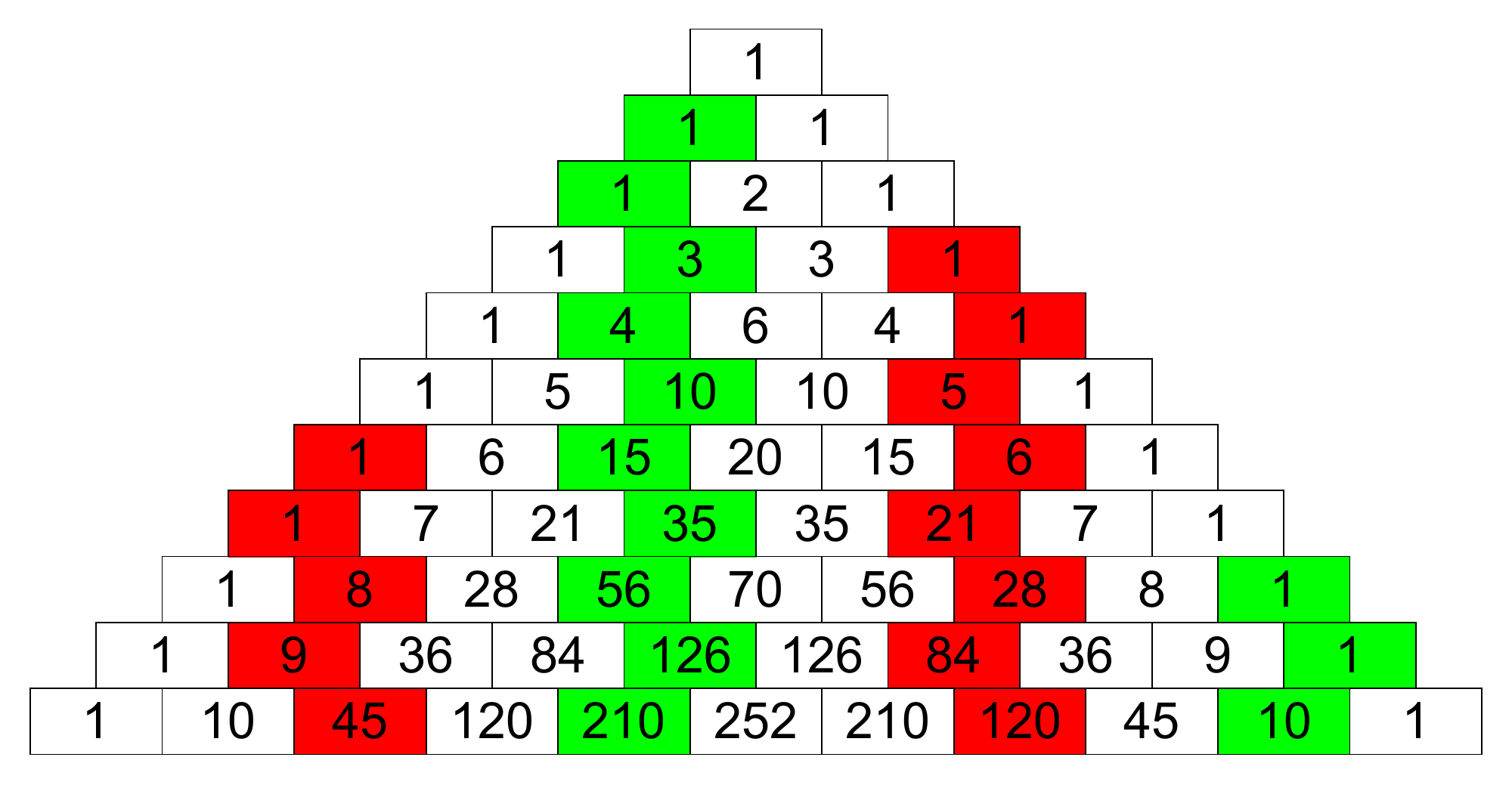}{fig:gs2}{Second GS identity.}

\begin{prop}[Second GS identity --- amateur version]
The sequence
\[
q(n)
=
\sum_{
-\floor{\frac{n-1}{5}}
\leq\lambda\leq
\floor{\frac{n+2}{5}}
}
(-1)^\lambda
\binom{n}{\floor{\frac{n-1+5\lambda}{2}}}
\]
satisfies the recurrence
\[
q(n) = q(n-1) + q(n-2)
,
\]
with
\[
q(0)=0;\;q(1)=1
.
\]
\end{prop}
The proof is the same.

Along with the GS identities comes the simpler `slalom identity':
Totalling the entries in any row of Figure \ref{fig:slalom} yields 1.
\fig{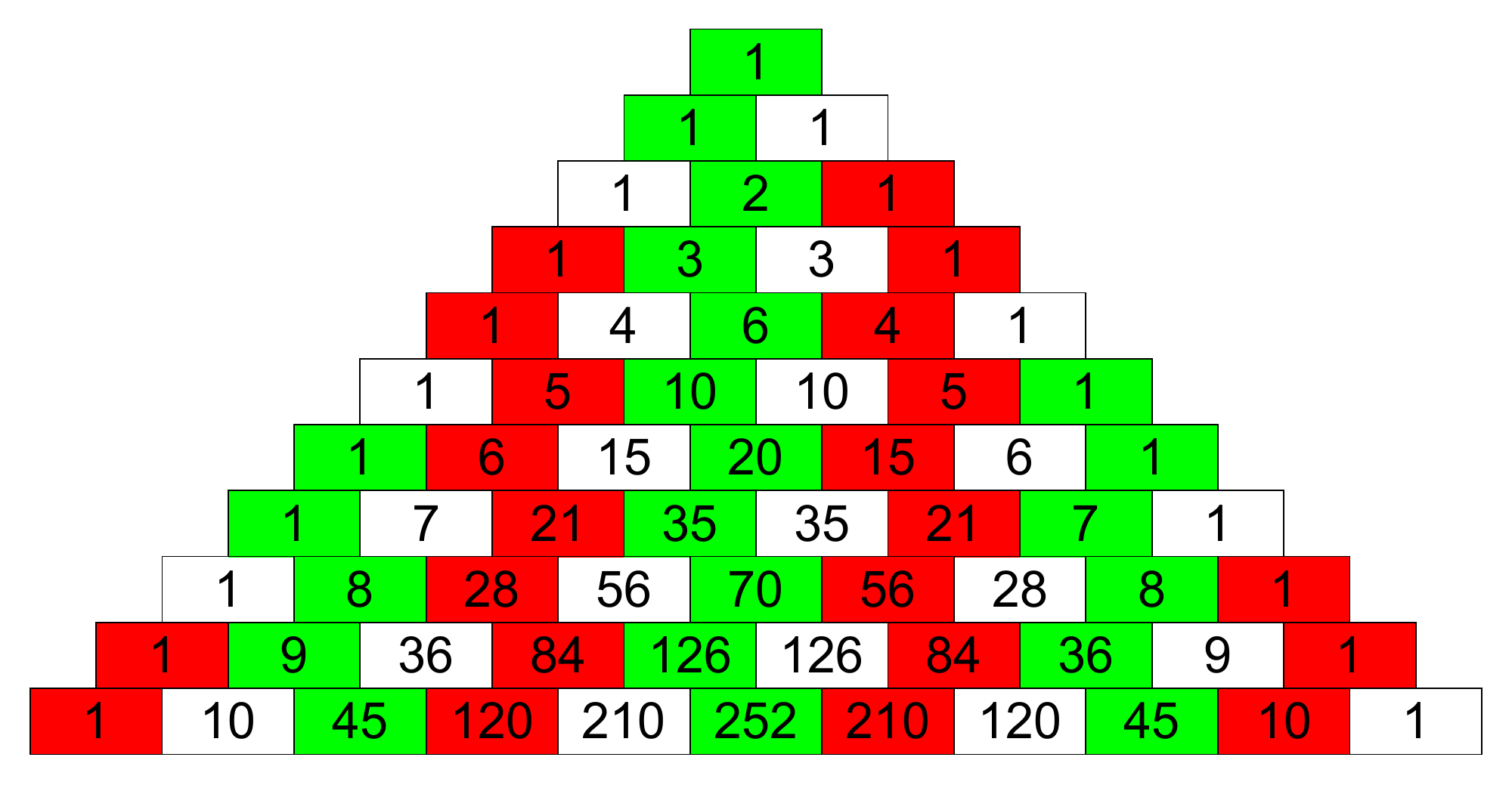}{fig:slalom}{Slalom identity.}
Here is the identity as Frobenius formulated it:
\begin{prop}[Slalom identity --- amateur version]
The sequence
\[
r(n)
=
\sum_{
-\floor{\frac{n}{3}}
\leq\lambda\leq
\floor{\frac{n+1}{3}}
}
(-1)^\lambda
\binom{n}{\floor{\frac{n+3\lambda}{2}}}
\]
is the constant sequence $1,1,1,\ldots$.
Or, as we prefer to say, it satisfies the recurrence
\[
r(n)=r(n-1)
,
\]
with
\[
r(0)=1
.
\]
\end{prop}

The proof of the slalom identity is like that of
the GS identities, only simpler:
See Figure \ref{fig:slalomproof}.
\fig{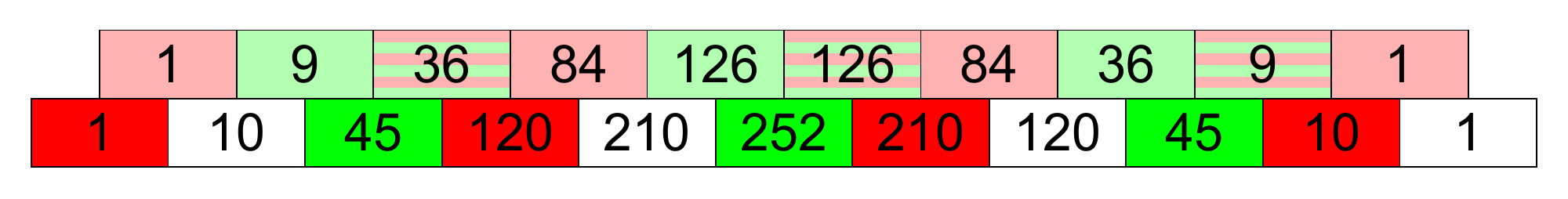}{fig:slalomproof}{Proof of the slalom identity.}
In symbols,
\begin{eqnarray*}
&&\binom{n}{k} - \binom{n}{k+1}
\\&=&
\binom{n-1}{k-1} + \binom{n-1}{k}
- \binom{n-1}{k} - \binom{n-1}{k+1}
\\&=&
\binom{n-1}{k-1}
 - \binom{n-1}{k+1}
.
\end{eqnarray*}

\section{Conclusion}

Shown the amateur alpine identities,
any number of combinatorialists will know to look for the
pro analogs involving Gaussian binomial coefficients,
and will easily find them;
the proofs precisely follow those of amateur versions
to which they reduce.
And then, just as Schur did, they will know how derive
the Rogers-Ramanujan identities from the GS identities,
either before or after deriving Euler's pentagonal number theorem
from the slalom identity.
We'll give details of this below, as Addenda.
We conclude here because, from the point of view of any number of 
combinatorialists, we're already done.

\section{Addendum: The pro alpine identities}

The polynomial analogs of binomial coefficients are the Gaussian
binomial coefficients, introduced by Gauss
\cite[p. 16]{gauss:binom}:
\[
\gbinom{n}{k}
=
\frac
{(x^n-1)(x^{n-1}-1)\ldots(x^{n-k+1}-1)}
{(x-1)(x^2-1)\ldots(x^k-1)}
\]
Here we are to understand that $\gbinom{n}{k}=0$ if any of $n,k,n-k$ are
negative.
(The notation $\gbinom{n}{k}$ was introduced by Schur
\cite[p. 128]{schur:rr}
in the
paper we're discussing.)

These are indeed polynomials because they satisfy the recurrence
\[
\gbinom{n}{k}
=
\gbinom{n-1}{k-1}+x^k \gbinom{n-1}{k}
\]
with
\[
\gbinom{0}{0}=1
.
\]
Because
\[
\gbinom{n}{k}=\gbinom{n}{n-k}
\]
we have the alternative recurrence
\begin{eqnarray*}
\gbinom{n}{k}
&=&
\gbinom{n}{n-k}
\\&=&
\gbinom{n-1}{n-k-1}
+
x^{n-k}
\gbinom{n-1}{n-k}
\\&=&
x^{n-k}
\gbinom{n-1}{k-1}
+
\gbinom{n-1}{k}
.
\end{eqnarray*}

\begin{prop}[First GS identity]
Setting 
\[
a(\lambda)
=
\frac{5\lambda^2-\lambda}{2},
\]
the sequence
\[
P(n)
=
\sum_{
-\floor{\frac{n}{5}}
\leq \lambda \leq
\floor{\frac{n+1}{5}}
}
(-1)^\lambda
x^{a(\lambda)}
\gbinom{n}{\floor{\frac{n+5\lambda}{2}}}
\]
satisfies the recurrence
\[
P(n) = P(n-1)+x^{n-1}P(n-2)
,
\]
with
\[
P(-1)=0;\;P(0)=1
.
\]
\end{prop}

\proofstart
We repeat the proof of the amateur identity,
sprinkling in appropriate powers of $x$ here and there.
\begin{eqnarray*}
\gbinom{n}{k-1}
&=&
\gbinom{n-1}{k-2}
+
x^{k-1}
\gbinom{n-1}{k-1}
\\&=&
\gbinom{n-1}{k-2}
+
x^{k-1}
\left(
x^{n-k}
\gbinom{n-2}{k-2}
+
\gbinom{n-2}{k-1}
\right)
\\&=&
\gbinom{n-1}{k-2}
+
x^{n-1}
\gbinom{n-2}{k-2}
+
x^{k-1}
\gbinom{n-2}{k-1}
,
\end{eqnarray*}
while
\begin{eqnarray*}
\gbinom{n}{k+1}
&=&
x^{n-k-1}
\gbinom{n-1}{k}
+
\gbinom{n-1}{k+1}
\\&=&
x^{n-k-1}
\left(
\gbinom{n-2}{k-1}
+
x^{k}
\gbinom{n-2}{k}
\right)
+
\gbinom{n-1}{k+1}
\\&=&
x^{n-k-1}
\gbinom{n-2}{k-1}
+
x^{n-1}
\gbinom{n-2}{k}
+
\gbinom{n-1}{k+1}
.
\end{eqnarray*}
Taking a linear combination with coefficients $\alpha,-\beta$,
we get
\begin{eqnarray*}
\alpha
\gbinom{n}{k-1}
-
\beta
\gbinom{n}{k+1}
&=&
\alpha
\gbinom{n-1}{k-2}
-
\beta
\gbinom{n-1}{k+1}
+
\\&&
x^{n-1}
\left(
\alpha
\gbinom{n-2}{k-2}
-
\beta
\gbinom{n-2}{k}
\right)
+
\\&&
(
\alpha
x^{k-1}
-
\beta
x^{n-k-1}
)
\gbinom{n-2}{k-1}
\
.
\end{eqnarray*}

Now consider any pair of entries $\gbinom{n}{k-1}$ and $\gbinom{n}{k+1}$ 
occupying cells of opposite color, meaning that
\[
k-1 = \floor{\frac{n+5\lambda}{2}}
;\;
k+1 = \floor{\frac{n+5(\lambda+1)}{2}}
.
\]
This happens just when $n$ and $\lambda$ have the same parity,
making
\[
k-1 = \frac{n+5\lambda}{2}
\]
and
\[
\lambda = \frac{2k-n-2}{5}
,
\]
so that
\[
a(\lambda+1)-a(\lambda)
=
5\lambda+2
=
2k-n
.
\]
This means that if we take
$\alpha=x^{a(\lambda)}$,
$\beta=x^{a(\lambda+1)}$
in the expressions above,
the final term vanishes because the coefficient is
\[
x^{a(\lambda)}
x^{k-1}
-x^{a(\lambda+1)}
x^{n-k-1}
=
x^{a(\lambda)}x^{k-1}
(1-x^{a(\lambda+1)-a(\lambda)}x^{n-2k})
= 
0.
\]
This establishes the recurrence.
\proofend

\begin{prop}[Second GS identity]
Setting 
\[
b(\lambda) = \frac{5\lambda^2-3\lambda}{2},
\]
the sequence
\[
Q(n)
=
\sum_{
-\floor{\frac{n-1}{5}}
\leq\lambda\leq
\floor{\frac{n+2}{5}}
}
(-1)^\lambda
x^{b(\lambda)}
\gbinom{n}{\floor{\frac{n-1+5\lambda}{2}}}
\]
satisfies the recurrence
\[
Q(n) = Q(n-1) + Q(n-2)
,
\]
with
\[
Q(0)=0;\;Q(1)=1
.
\]
\end{prop}

\proofstart
Same as above, only now
the entries $\gbinom{n}{k-1}$ and $\gbinom{n}{k+1}$ have opposite color
just when $n$ and $\lambda$ have opposite parity,
making
\[
k-1 = \frac{n-1+5\lambda}{2}
\]
and
\[
\lambda = \frac{2k-1-n}{5}
,
\]
so that once again
\[
b(\lambda+1)-b(\lambda) = 5 \lambda +1 = 2k-n
.
\mathproofend
\]

\begin{prop}[Slalom identity]
Setting
\[
c(\lambda) =
\frac{3\lambda^2-\lambda}{2}
,
\]
the sequence
\[
R(n)
=
\sum_{
-\floor{\frac{n}{5}}
\leq \lambda \leq
\floor{\frac{n+1}{5}}
}
(-1)^\lambda
x^{c(\lambda)}
\gbinom{n}{\floor{\frac{n+3\lambda}{2}}}
\]
is the constant sequence $1,1,1,\ldots$.
Or, as we prefer to say, it satisfies the recurrence
\[
R(n)=R(n-1)
,
\]
with
\[
R(0)=1
.
\]
\end{prop}

\proofstart
\begin{eqnarray*}
&&\alpha \gbinom{n}{k} - \beta \gbinom{n}{k+1}
\\&=&
\alpha
\left(
\gbinom{n-1}{k-1}
+
x^k
\gbinom{n-1}{k}
\right)
-
\beta
\left(
x^{n-k-1}
\gbinom{n-1}{k}
+
\gbinom{n-1}{k+1}
\right)
\\&=&
\alpha \gbinom{n-1}{k-1}
-
\beta
\gbinom{n-1}{k+1}
+
\\&&
(\alpha x^k - \beta x^{n-k-1})
\gbinom{n-1}{k}
\end{eqnarray*}
If the entries $\gbinom{n}{k}$ and $\gbinom{n}{k+1}$
are in cells of opposite color,
then $n$ and $\lambda$ have the same parity, making
\[
k=
\frac{n+3\lambda}{2}
\]
and
\[
\lambda = \frac{2k-n}{3}
,
\]
so that
\[
c(\lambda+1)-c(\lambda)
=
3\lambda+1
=
2k-n+1,
\]
so that with $\alpha=c(\lambda)$, $\beta=c(\lambda+1)$,
the coefficient of $\gbinom{n-1}{k}$ above is
\[
x^{c(\lambda)} x^k - x^{c(\lambda+1)}x^{n-k-1}
=
x^{c(\lambda)} x^k (1 - x^{c(\lambda+1)-c(\lambda)}x^{n-2k-1})
=
0
.
\mathproofend
\]

\section{Addendum:  Euler and Rogers-Ramanujan}

Schur deduces the Euler pentagonal number theorem and
Rogers-Ramanujan identities by combining the alpine identities
with the Jacobi triple product formula
(cf.\ Gauss
\cite{gauss:100})
in a way that by now
is very well known
(cf. Andrews and Eriksson \cite{andrewseriksson}).
We briefly review this here.

Attribute to an integer partition $a = (a_1,\ldots,a_n)$,
\[
a_1 \geq a_2 \geq \ldots \geq a_n \geq 1
\]
the \emph{weight} $x^{a_1+\ldots+a_n}$.
Totaling the weights of all partitions gives
\[
E = \prod_i \frac{1}{(1-x^i)}
.
\]
The Gaussian binomial $\gbinom{n}{k}$ totals the weights
of partitions with at most $n-k$ parts, each of size at most $k$.
If $k,n-k$ both tend to infinity, we get $E$:
\[
\lim_{k,n-k\goesto \infty} \gbinom{n}{k}
=
E
.
\]

Taking $n$ to infinity in the slalom identity,
we get
\[
1 =
\sum_\lambda
(-1)^\lambda
x^{\frac{3\lambda^2-\lambda}{2}}
E
,
\]
or
\[
\prod_i (1-x^i) =
\sum_\lambda
(-1)^\lambda
x^{\frac{3\lambda^2-\lambda}{2}}
.
\]
This is Euler's pentagonal number theorem.

Call $a$ 
a \emph{kangaroo partition} if
\[
\min(a_1-a_2,a_2-a_3,\ldots,a_{n-1}-a_n)
\geq 2
.
\]
$P(n)$ is the weight sum
for kangaroo partitions with maximum part at most $n-1$.
Taking $n$ to infinity in GS1, on the left we get
the weight sum $P(\infty)$ over
all kangaroo partitions,
while on the right all the Gaussian binomials turn into $E$:
\[
P(\infty) 
=
J E
,
\]
where
\[
J =
\sum_\lambda
(-1)^\lambda
x^{\frac{5\lambda^2-\lambda}{2}}
.
\]
But from the Jacobi triple product identity we have
\[
J =
\prod_{k \geq 0} (1-x^{5k+2})(1-x^{5k+3})(1-x^{5k+5})
,
\]
so
\[
P(\infty)
=
\prod_{k \geq 0} \frac{1}{(1-x^{5k+1})(1-x^{5k+4})}
.
\]
The right hand side here enumerates partitions with all parts congruent to
$\pm1$ mod $5$.
This is the first Rogers-Ramanujan identity.

We get the second Rogers-Ramanujan identity from GS2 in like manner:
\[
Q(\infty)
=
\prod_{k \geq 0} (1-x^{5k+2})(1-x^{5k+3})
.
\]
On the left we have kangaroo partitions with minimum part size at least $2$,
while on the right we have partitions into parts congruent to $\pm2$ mod $5$.

\section{Addendum:  Bijective proofs of the alpine identities}

From our algebraic proofs of the alpine identities we can extract
bijections that pair up terms of opposite sign, leaving 
a single positive term in each row for the slalom, or a Fibonacci number
of positive terms
for GS1 and GS2.

For GS1, if you trace
through what cancels with what,
you will find that
the terms that remain correspond to partitions built from L-shaped
pieces with equal or all-but-equal prongs, as indicated in Figure
\ref{fig:fib}.
All other partitions are matched up into pairs of opposite sign
as shown in Figure
\ref{fig:cancel}.
Here's pseudo-code for bijections that do this pairing:
\begin{verbatim}
leftmatch[{A, steps___}] := Join[rightmatch[{steps}], {B}]
leftmatch[{B, steps___, A}] := Join[{B}, leftmatch[{steps}], {A}]
leftmatch[{B, steps___, B}] := {A, steps, A}
leftmatch[{steps___}] := {steps}

rightmatch[{steps___, B}] := Join[{A}, leftmatch[{steps}]]
rightmatch[{B, steps___, A}] := Join[{B}, rightmatch[{steps}], {A}]
rightmatch[{A, steps___, A}] := {B, steps, B}
rightmatch[{steps___}] := {steps}
\end{verbatim}
Apply {\tt rightmatch} or {\tt leftmatch} according as the parity of
$n + \lambda$ is even or odd.

\fig{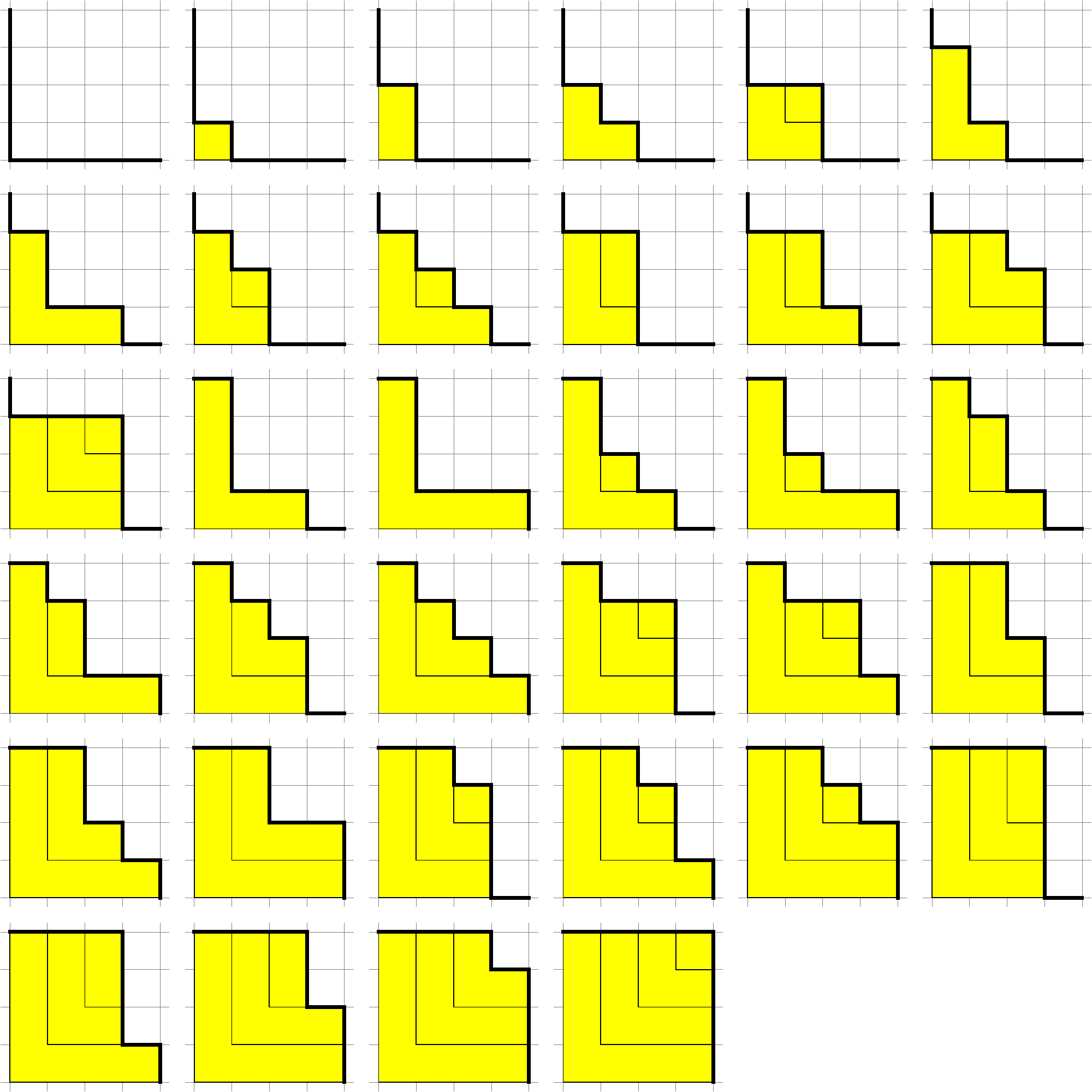}{fig:fib}{The 34 integer partitions belonging to $P(8)$.}

\figsize{390pt}{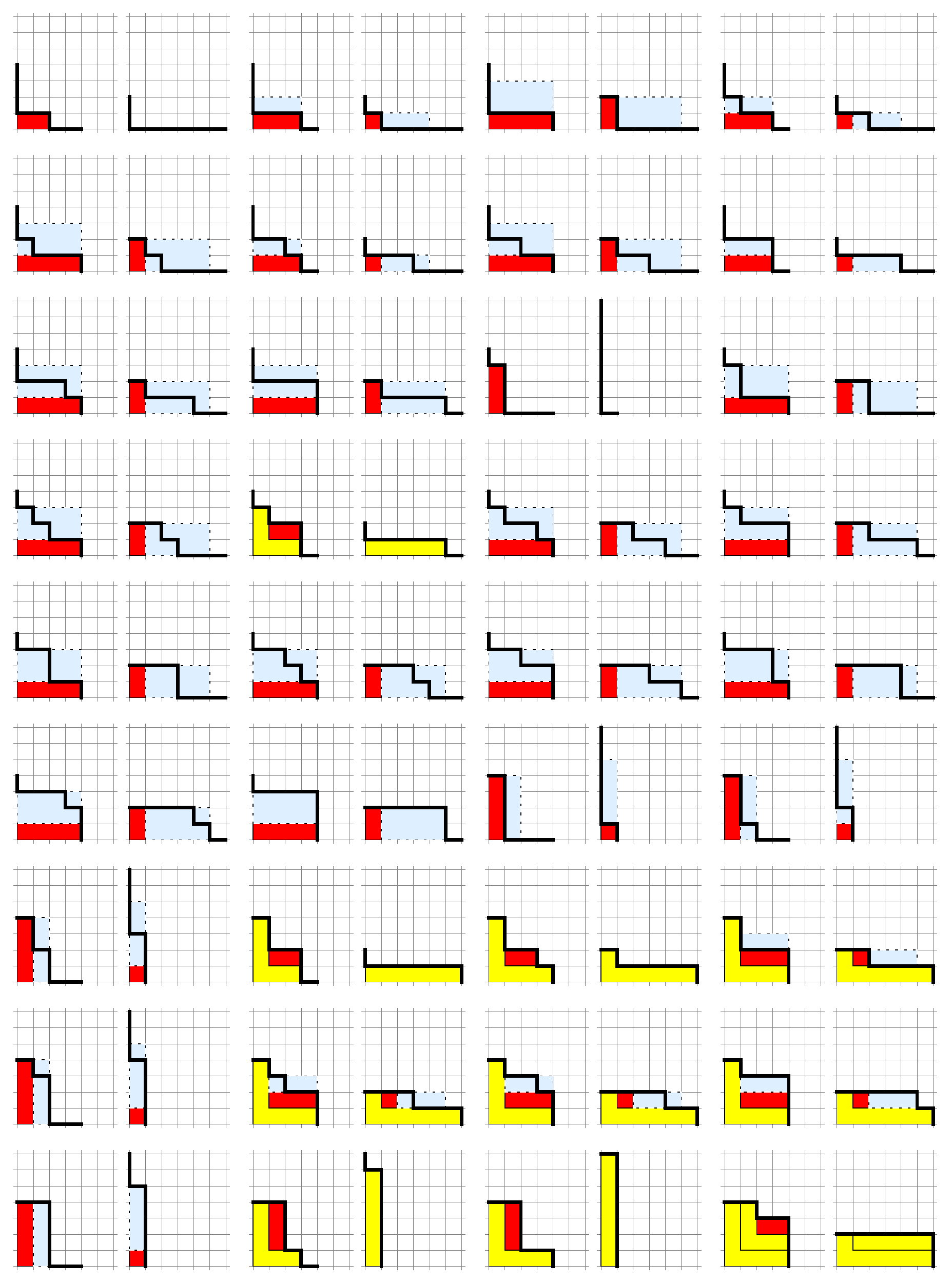}{fig:cancel}{The 36 canceling pairs on the right of the
GS1 identity for P(8).}

Verifying directly that these bijections have the stated properties
gives us bijective proofs of the alpine identities.
It would be a mistake to regard these bijective proofs
as distinct from the algebraic proofs
from which they are extracted.
They are, collectively, just another manifestation of
Frobenius's last proof.

\pagebreak
\bibliography{flp}
\bibliographystyle{hplain}

\end{document}